\documentclass[a4paper,12pt]{amsart}
\usepackage{amsmath,amscd}
\newtheorem{thm}{Theorem}[section]

\newtheorem{lem}[thm]{Lemma}
\newtheorem{cor}[thm]{Corollary}
\newtheorem{prop}[thm]{Proposition}

\theoremstyle{definition}

\newtheorem{setup}[thm]{}
\newtheorem{exmp}[thm]{Example}

\newtheorem{rem}[thm]{Remark}
\theoremstyle{remark}

\title{Weak boundedness theorems for canonically fibered Gorenstein minimal
3-folds}
\author{Meng Chen}
\address{Institute of Mathematics, Fudan University, Shanghai, 200433, PR China}
\email{mchen@fudan.edu.cn}
\date{}
\thanks{This paper is
supported by the National Natural Science Foundation of China
(No.10131010), Shanghai Scientific $\&$ Technical Commission
(Grant 01QA14042) and SRF for ROCS, SEM}
\begin{document}

\begin{abstract}
Let $X$ be a Gorenstein minimal projective 3-fold with at worst locally factorial terminal
singularities. Suppose the canonical map is of fiber type. Denote by $F$ a smooth model of
a generic irreducible element in fibers of $\phi_1$ and so $F$ is a curve or a smooth surface.
The main result is that there is a computable constant $K$ independent of $X$ such that
$g(F)\le 647$
or $p_g(F)\le 38$ whenever $p_g(X)\ge K$.
\end{abstract}
\maketitle

\pagestyle{myheadings}
\markboth{\hfill M. Chen\hfill}{\hfill Boundedness Theorems\hfill}

\section*{\rm Introduction}

In this paper, we study the boundedness problem of a canonical family of curves or
surfaces on certain algebraic 3-folds. Let $X$ be a Gorenstein minimal projective
3-folds with locally factorial terminal singularities. Suppose $p_g(X)\ge 2$. One may define
the canonical map $\phi_1$. Denote $\alpha:=\dim\phi_1(X)$. We say that $\phi_1$ is {
\it of fiber type} if $\alpha<3$. Denote by $F$ a smooth model of a generic irreducible
element in fibers of
$\phi_1$. It is interesting to know whether the birational invariants of $F$ are
bounded.
For technical reason, we have to assume on $X$ an effective Miyaoka-Yau inequality. This is
the reason that we
are only able to treat a Gorenstein minimal 3-folds by virtue of (\cite{Mi}). Apparently
we have $\kappa(X)\ge 1$. If $\kappa(X)\le 2$, by the Base Point Free Theorem (\cite{KMM, K-M}),
one can see
that $F$ is special. So, throughout, we may suppose $\kappa(X)=3$, i.e. $X$ is of general
type. What we can do is as follows (see Theorem 1.5, Theorem 2.3 and Theorem 2.8).

\begin{thm}\label{T:0.1} There is a computable constant $K$. For all $X$, a Gorenstein
minimal projective 3-fold with locally factorial terminal singularities, denote by $\phi_1$ the canonical map
of $X$ and by $F$ a smooth model of a generic irreducible element in fibers of $\phi_1$.
The following holds:

(1) If $\dim\phi_1(X)=2$, then $F$ is a curve with geometric genus $g(F)\le 647$
provided $p_g(X)\ge K$.

(2) If $\dim\phi_1(X)=1$, then  $F$ is a surface of general type with geometric genus
$p_g(F)\le 38$ provided $p_g(X)\ge K$.
\end{thm}

\begin{rem}In the above theorem, we don't know the boundedness when $p_g(X)<K$.
The relations among $p_g(X)$, $q(X)$ and $h^2({\mathcal O}_{X})$ are expected to know to solve the
problem. We also believe
that those bounds in \ref{T:0.1} might be far from sharp.
\end{rem}

I thank a referee for effective suggestions. This paper was
finally revised while I was visiting Kang Zuo at the Chinese
University of Hong Kong. Special thanks are due to Kang Zuo for
many helps and the hospitality.

\section{\rm Canonical family of surfaces}\label{Sec:1}

\begin{setup}\label{S:1.1}  {\bf Setting up.}
Let $X$ be a minimal projective 3-fold with ${\Bbb Q}$-factorial terminal singularities.
Suppose $p_g(X)\ge 2$. We can define the canonical map $\phi_1$ which is usually a rational map.
Take the birational modification $\pi: X'\longrightarrow X$, according to Hironaka, such that

(i) $X'$ is smooth;

(ii) the movable part of $|K_{X'}|$ is base point free;

(iii) certain $\pi^*(K_X)$ has supports with only normal crossings.

Denote by $g$ the composition $\phi_1\circ\pi$. So
$g: X'\longrightarrow W'\subseteq{\Bbb P}^{p_g(X)-1}$
is a morphism. Let
$g: X'\overset{f}\longrightarrow W\overset{s}\longrightarrow W'$
be the Stein factorization of $g$. We can write
$$K_{X'}=\pi^*(K_X)+E=S_1+Z_1,$$
where $S_1$ is the movable part of $|K_{X'}|$, $Z_1$ the fixed part and $E$ is an effective
${\Bbb Q}$-divisor which is a sum of distinct exceptional divisors.

If $\dim\phi_1(X)<3$, $f$ is a called a {\it derived fibration of}
$\phi_1$. We mean {\it a generic irreducible element in fibers of} $\phi_1$ by
a general fiber of $f$.

If $\dim\phi_1(X)=2$, a general fiber $C$ of $f$ is a smooth curve of genus $g(F)$ which should not be
confused with the morphism $g$.

If $\dim\phi_1(X)=1$, a general fiber $F$ of $f$ is a smooth projective surface of
general type. Denote by $F_0$ the minimal model of $F$. Denote by $b$ the genus of the smooth
curve $W$.
\end{setup}

\begin{prop}\label{P:1.2} Keep the above notations. Suppose $X$ is Gorenstein minimal of
general type with locally factorial terminal singularities,
$\dim\phi_1(X)=1$ and $b=0$.
If $p_g(X)\ge 2k+2$ for certain $k\ge 4$, then there is an effective ${\Bbb Q}$-divisor $E_k'$
on $X'$
such that
$$\pi^*(K_X)|_F=_{\Bbb Q}\frac{k}{k+1}\sigma^*(K_{F_0})+E_k'$$
where $\sigma: F\rightarrow F_0$ is contraction onto the minimal model.
\end{prop}
\begin{proof}
Let $M_{k+1}$ be the movable part of $|(k+1)K_{X'}|$. Then we can write
$$(k+1)\pi^*(K_X)\sim M_{k+1}+E_{k+1}'$$
where $E_{k+1}'$ is an effective divisor with respect to $k$.
Therefore we see that
$(k+1)\pi^*(K_X)|_F\ge M_{k+1}|_F$. Let $N_k$ be the movable part of $|kK_F|$.
According to \cite{Bo}, $|kK_{F_0}|$ is base point free.  Thus
$N_k=\sigma^*(kK_{F_0})$ where $\sigma:F\longrightarrow F_0$ is contraction onto the minimal model.

We use the approach in \cite{Kol} (Corollary 4.8)  to prove $M_{k+1}|_F\ge N_k$. In fact, we have a fibration
$f:X'\longrightarrow {\Bbb P}^1$. Because
$p_g(X)\ge 2k+2$,
we can see that ${\mathcal O}(2k+1)\hookrightarrow f_*\omega_{X'}$. Thus we have
$${\mathcal  E}:={\mathcal O}(1)\otimes f_*\omega_{X'/{\Bbb P}^1}^k=
{\mathcal O}(2k+1)\otimes f_*\omega_{X'}^k\hookrightarrow f_*\omega_{X'}^{k+1}.$$
Note that
$H^0({\Bbb P}^1, f_*\omega_{X'}^{k+1})\cong H^0(X', \omega_{X'}^{k+1}).$
It is well known that ${\mathcal E}$ is generated by global sections and that
$f_*\omega_{X'/{\Bbb P}^1}^k$ is a sum of line bundles with non-negative degree
(\cite{F, V2, V3}). Thus the global sections of ${\mathcal  E}$ can
distinguish
different fibers of $f$. On the other hand, the local sections of $f_*\omega_{X'}^k$ give
the k-canonical map of $F$ and these local sections can be extended to global sections of
${\mathcal  E}$. This means $M_{k+1}|_F\ge N_k$. The proposition is proved.
\end{proof}

\begin{thm}\label{T:1.3} \text{\bf (\cite{Ch1}, Theorem 1)} Let $X$ be a minimal projective
3-fold of general type with at worst ${\Bbb Q}$-factorial terminal singularities.
Suppose $\dim\phi_1(X)=1$. If $b\ge 2$, then either
$$p_g(F)=1\ \text{and}\  \ p_g(X)\ge b-1\ \ \text{or}$$
$$b=p_g(F)=p_g(X)=2.$$
\end{thm}

\begin{setup}\label{S:1.4} {\bf Reduction to the case $b=0$}
We mainly study a Gorenstein minimal 3-fold.
By virtue of \ref{T:1.3}, one only has to consider the case $b\le 1$. If $b=1$, then
$p_g(F)\le 38$ according to \cite{Ch1} (Theorem 2(1)). {}From now on within this section, we
may suppose $b=0$. So $W={\Bbb P}^1$.
\end{setup}

\begin{thm}\label{T:1.5} There exists a constant $K_1\le 800$ such that, for all
$X$, a Gorenstein minimal projective 3-fold of general type with at worst locally factorial
terminal singularities, satisfying $\dim\phi_1(X)=1$, the following holds.

(1) $p_g(F)\le 38$ provided $p_g(X)\ge K_1$.

(2) $p_g(F)\le 5002$ provided $p_g(X)\ge 23$.
\end{thm}
\begin{proof} We only have to prove the theorem for the case $b=0$. In this case, $W={\Bbb P}^1$. We have a derived fibration $f: X'\longrightarrow {\Bbb P}^1$. We have
$$q(X)=h^1({\Bbb P}^1, R^1f_*\omega_{X'})\le q(F).$$
Because $F$ is of general type, we have
$p_g(F)\ge 2q(F)-4$
according to \cite{Be2}. Thus the Miyaoka-Yau inequality (\cite{Mi}) becomes
\begin{align*}
K_X^3&\le -72[1-q(X)+h^2({\mathcal O}_{X})-p_g(X)]\\
&\le 72[(p_g(X)-1)+q(F)]\\
&\le 72[p_g(X)-1]+36p_g(F)+144.
\end{align*}
When $p_g(X)\ge 2k+2$ for $k\ge 4$, Proposition 1.2 gives
$$K_X^3\ge \frac{k^2}{(k+1)^2}\cdot K_{F_0}^2\cdot [p_g(X)-1].$$
Denote $A:=p_g(X)-1$ and
$B:=\frac{k^2}{(k+1)^2}.$
Also using the inequality $K_{F_0}^2\ge 2p_g(F)-4$, we obtain the following inequality
$$
p_g(F)\le\frac{36A+2AB+72}{AB-18}
=\frac{(36+2B)+\frac{72}{A}}{B-\frac{18}{A}}.$$
Noting that $A\rightarrow +\infty$ and $B\rightarrow 1$ if
$k\rightarrow +\infty$, we can easily find a constant $K_1$ such that
$p_g(F)\le 38$ whenever
$p_g(X)\ge K_1$. The proof is completed.
\end{proof}

\begin{rem}\label{R:1.6}
Among known examples, the biggest possible $p_g(F)$ is 5 (\cite{C-C}, Example 2(e))
where the total space is a smooth minimal 3-fold and the canonical system is composed of a
pencil of surfaces over rational curve. Can one construct more examples with bigger $p_g(F)$?
\end{rem}

\section{\rm Canonical family of curves}\label{S:2}

\begin{setup}\label{S:2.1} {\bf Setting up.}
Within this section, we keep the same notations as in \ref{S:1.1} and suppose $X$ is a Gorenstein minimal
3-fold of general type with at worst locally factorial terminal singularities and
$\dim\phi_1(X)=2$.
In this case, we may even suppose $W$ is a smooth surface and so a general fiber $C$ of
$f$ is a smooth projective curve.
\end{setup}

\begin{setup}\label{S:2.2} {\bf An inequality.}
We need a nontrivial inequality in terms of $g(C)$ as
$$K_X^3\ge \frac{2}{3}(g(C)-1)(p_g(X)-2)$$
according to \cite{Ch3} (Theorem 4.1(ii)). We omit the proof here simply
because, on one hand side, this type of inequality is expectable and, on the other hand,
the proof is quite clear and short there.
\end{setup}

\begin{thm}\label{T:2.3}  There is a computable constant $K_2$ such that for all $X$, a
Gorenstein minimal projective 3-fold of general type with at worst locally factorial
terminal singularities, satisfying $\dim\phi_1(X)=2$ and $\kappa(W)\ge 0$,
the following holds:

(1) $g(C)\le 109$ provided $p_g(X)\ge K_2$.

(ii) $g(C)\le 11989$ provided $p_g(X)>109$.
\end{thm}
\begin{proof}
We have a canonically derived fibration $f: X'\longrightarrow W$.  According
to \cite{Kol} and
\cite{V2, V3}, we know that
$R^1f_*\omega_{X'}\cong \omega_{W}$ and $f_*\omega_{X'/W}$ is torsion free.
We have
\begin{align*}
q(X)&=h^1(W, R^1f_*\omega_{X'})+h^2(W, f_*\omega_{X'})\\
&=q(W)+h^2(W, f_*\omega_{X'})\le q(W)+g(C).
\end{align*}
\begin{align*}
h^2({\mathcal O}_{X})&=h^0(W, R^1f_*\omega_{X'})+h^1(W, f_*\omega_{X'})\\
&=p_g(W)+h^1(W, f_*\omega_{X'}).
\end{align*}
On the surface $W$, we have
$q(W)-p_g(W)=1-\chi({\mathcal O}_{W}).$
By virtue of the birational classification theory on surfaces, we can see that
$$\chi({\mathcal O}_{W})\ge 0\ \ \text{whenever}\ \  \kappa(W)\ge 0.$$
Thus we definitely have
$q(W)-p_g(W)\le 1$
under the assumption of the theorem. So we get
$$q(X)-h^2({\mathcal O}_{X})\le q(W)-p_g(W)+g(C)\le g(C)+1.$$
$$\chi(\omega_X)=p_g(X)+q(X)-h^2({\mathcal O}_{X})-1\le p_g(X)+g(C).$$
By \ref{S:2.2} and Miyaoka's inequality, we have
$$\frac{2}{3}(g(C)-1)[p_g(X)-2]\le 72[p_g(X)+g(C)].$$
Denote $A:=p_g(X)-1$. The above inequality becomes
$$
g(C)\le\frac{109A+108}{A-108}=\frac{109+\frac{108}{A}}{1-\frac{108}{A}}.$$
It's easy to find a constant $K_2$ such that $g(C)\le 109$ provided $p_g(X)\ge K_2$.
\end{proof}

\begin{setup}\label{S:2.4} {\bf The case $\kappa(W)=-\infty$.}
We know that $W$ is birationally equivalent to a ruled surface
or ${\Bbb P}^2$. Denote $b_0:=q(W)$. This situation is more complicated. We have to do
some preparation in order to prove the boundedness.
\end{setup}

\begin{cor}\label{C:2.5} There exists a computable constant $K_3$ such that for all $X$,
 a Gorenstein minimal projective 3-fold of general type with at worst  locally factorial
terminal singularities, satisfying $\dim\phi_1(X)=2$ and $\kappa(W)=-\infty$, we have
$g(C)\le b_0+109$ provided $p_g(X)\ge K_3$.
\end{cor}
\begin{proof}
According to the proof of \ref{T:2.3}, we only have to estimate $\chi(\omega_X)$.
Because $q(X)\le b_0+g(C)$, we have
$$\chi(\omega_X)=p_g(X)+q(X)-h^2({\mathcal O}_{X})-1\le p_g(X)+b_0+g(C)-1.$$
By \ref{S:2.2} and Miyaoka's inequality, we have
$$\frac{2}{3}(g(C)-1)[p_g(X)-2]\le 72 [p_g(X)+b_0+g(C)-1].$$
Denote $A_1=p_g(X)-2$. The inequality is
$$g(C)\le 109\cdot \frac{A_1}{A_1-108}+b_0\cdot \frac{108}{A_1-108}+\frac{108}{A_1-108}.$$
Noting that
$$\frac{A_1}{A_1-108}\mapsto 1, \ \frac{108}{A_1-108}\mapsto 0\ \ \text{when}\ \ p_g(X)\mapsto +\infty,$$
we can easily find a constant $K_3$ such that $g(C)\le b_0+109$ provided $p_g(X)\ge K_3$.
\end{proof}

\begin{lem}\label{L:2.6} Let $X$ be a Gorenstein minimal projective 3-fold of general type
with at worst locally factorial terminal singularities.
Suppose $\dim\phi_1(X)=2$, $\kappa(W)=-\infty$ and $b_0\ge 3$. Then
$p_g(S_1)\le \frac{9}{2}K_X^3-2p_g(X)+6$ where $S_1$ is by definition the movable part of
$|K_{X'}|$.
\end{lem}
\begin{proof}
Because $W$ is birationally ruled, we have a projection map $p:W\longrightarrow B$ onto a smooth
curve $B$ with genus $b_0$. Because $b_0>0$, $p$ is a morphism and a fibration. Because
$\phi_1$ factors through $f$ and $\dim\phi_1(X)=2$,
$\dim\ f(S_1)=1$. This means $f$ maps a general $S_1$ onto a curve $H$ in $W$ and $H$ is exactly the pull back of a hyperplane section $H'$ in $W'$ under the map $s:W\longrightarrow W'$. We see that $H$ is a nef and big divisor on $W$ and $H$ is a smooth ireducible curve. Therefore there is a finite map $\theta: H\longrightarrow B$, because $p$ doesn't contract $H$. So $b_1:=g(H)\ge g(B)=b_0$. We have the following commutative diagram:
$$\begin{CD}
S_1 @>{\alpha=f|_{S_1}}>> H @>\theta>>  B\\
@V{j_1}VV   @VV{j_2}V    @VV\|V\\
X' @>>f> W  @>>\rho> B\\
@V{\pi}VV     @VV{s}V     &\\
X  @>>{\text{the rational map}\ \ \phi_1}>{W'\subset{\Bbb P}^N}
\end{CD}$$
where $j_1$ and $j_2$ are both inclusion. We can see that $\alpha: S_1\longrightarrow H$ is a
fibration where a general fiber $C'$ lies in the same numerical class as a fiber of $f$.
Thus $g(C')=g(C)$. Because $b_0\ge 3$, $q(S_1)\ge b_1\ge b_0\ge 3$. So $|K_{S_1}|$ defines
a generically finite map onto its image according to \cite{X}. Denote by $N_1$ the movable part of $|K_{S_1}|$. Then, according to Noether's inequality,
$N_1^2\ge 2p_g(S_1)-4.$
We can bound $p_g(S_1)$ by studying $\phi_{3}$. It's obvious that
$|K_{X'}+\pi^*(K_X)+S_1|\subset |3K_{X'}|.$
Denote by $M_3'$ the movable part of $|K_{X'}+\pi^*(K_X)+S_1|$ and by $M_3$ the movable part of $|3K_{X'}|$. Then we have
$3\pi^*(K_X)\ge M_3\ge M_3'.$
The vanishing theorem gives
$$|K_{X'}+\pi^*(K_X)+S_1|\bigm|_{S_1}=|K_{S_1}+L|\supset|N_1+L_0|,$$
where $L=\pi^*(K_X)|_{S_1}$ and $L_0=S_1|_{S_1}$.
Because $S_1^3=0$, we can see that $|L_0|$ is composed of a free pencil of curves on the
surface $S_1$. We can write $L_0\equiv a_1C'$ where $a_1\ge p_g(X)-2$. Since $N_1+L_0$ is
movable, we get from \cite{Ch2} (Lemma 2.7) that
$$3L\ge M_3'|_{S_1}\ge N_1+L_0.$$
Thus $9L^2\ge (N_1+L_0)^2$. Noting that $N_1\cdot C'\ge 2$, we can easily derive the inequality
$$p_g(S_1)\le\frac{9}{2}K_X^3-2p_g(X)+6.$$
\end{proof}

\begin{prop}\label{P:2.7}  Let $X$ be a Gorenstein minimal projective 3-fold of general
type with at worst locally factorial terminal singularities.
Suppose $\dim\phi_1(X)=2$, $\kappa(W)=-\infty$ and $b_0\ge 3$. If $g(C)\ge 648$, then
$b_0\le p_g(X)+653.$
\end{prop}
\begin{proof} According to the proof of Lemma \ref{L:2.6},
we have a surface fibration $\alpha: S_1\longrightarrow H$ where $g(C')=g(C)$ and $b_1=g(H)\ge b_0$.
We shall use a filtration of vector bundles on $H$ which was first studied
by G. Xiao (\cite{X}).
Considering the natural map
$$H^0(X', S_1)\overset{\gamma_1}\longrightarrow \Lambda_1\subset H^0(S_1, L_0)$$
where $\Lambda_1$ is the image of $\gamma_1$,
we can see that $\alpha$ is exactly obtained by taking the Stein factorization of
$\Phi_{\Lambda_1, S_1}$. Note that $L_0\le K_{S_1}$ and
$H^0(S_1, \omega_{S_1})\cong H^0(H,\alpha_*\omega_{S_1})$ where $\alpha_*\omega_{S_1}$ is
a vector bundle of rank $g(C)$ on the curve $H$.
Denote by ${\mathcal  L}_0$ the saturated sub-bundle of $\alpha_*\omega_{S_1}$, which is generated by sections in $\Lambda_1$. Then it's obvious that ${{\mathcal O}_ L}_0$ is a line bundle because $|L_0|$ is composed of fibers of $\alpha$. So we get an extension of $\alpha_*\omega_{S_1}$:
$$0\longrightarrow {\mathcal  L}_0\longrightarrow
\alpha_*\omega_{S_1}\longrightarrow {\mathcal  L}_1\longrightarrow 0.$$
Because $\alpha_*\omega_{S_1/H}$ is semi-positive, we have
$\deg({\mathcal L}_1\otimes \omega_H^{-1})\ge 0$
 i.e. $\deg({\mathcal  L}_1)\ge 2(g(C)-1)(b_1-1).$
The R-R gives
$h^0(H, {\mathcal  L}_1)\ge (g(C)-1)(b_1-1).$
Noting that $\deg({\mathcal  L}_0)>0$, we get, by applying the Clifford's theorem,
that
$h^1(H, {\mathcal L}_0)\le b_1-1.$
{}From the long exact sequence
$$0\longrightarrow H^0({\mathcal  L}_0)\longrightarrow
H^0(\alpha_*\omega_{S_1})\longrightarrow H^0({\mathcal  L}_1)\longrightarrow
H^1({\mathcal L}_0)\longrightarrow\cdots,$$
we have the inequality
$$h^0({\mathcal  L}_0)-h^0(\alpha_*\omega_{S_1})+h^0({\mathcal  L}_1)-
h^1({\mathcal  L}_0)\le 0.$$
Noting that $h^0({\mathcal  L}_0)\ge 2$ and using Lemma \ref{L:2.6}, we have
$$2+(g(C)-2)(b_1-1)\le p_g(S_1)\le \frac{9}{2}K_X^3-2p_g(X)+6.$$
Also noting that $b_1\ge b_0$ and using Miyaoka's inequality, one can get
$$(g(C)-2)(b_0-1)\le 322p_g(X)+324b_0+324g-320,$$
$$(g(C)-326)b_0\le 322p_g(X)+325g-322.$$
If $g(C)>326$, we have
$$b_0\le \frac{322p_g(X)}{g(C)-326}+325+\frac{105628}{g(C)-326}.$$
It's easy to see $b_0\le p_g(X)+653$ provided $g(C)\ge 648$.
\end{proof}

\begin{thm}\label{T:2.8}  There exists a computable constant $K_4$ such that for all
$X$, a Gorenstein minimal projective 3-fold of general type with at worst locally
factorial terminal singularities, satisfying $\dim\phi_1(X)=2$ and $\kappa(W)=-\infty$,
we have $g(C)\le 647$ provided $p_g(X)\ge K_4$.
\end{thm}
\begin{proof}
If $b_0\le 2$, then Corollary \ref{C:2.5} gives $g(C)\le 111$ provided
$p_g(X)\ge K_3$. If $b_0>2$
and $g(C)\ge 648$, we have $b_0\le p_g(X)+653$ by Proposition \ref{P:2.7}.
Applying \ref{S:2.2} and
Miyaoka-Yau inequality, we obtain
\begin{align*}
\frac{2}{3}(g(C)-1)[p_g(X)-2]&\le 72[p_g(X)+b_0+g(C)-1]\\
&\le 144p_g(X)+72g(C)+72\cdot 652.
\end{align*}
If $p_g(X)>110$, we can obtain
$$g(C)\le 217\cdot\frac{p_g(X)}{p_g(X)-110}+\frac{70414}{p_g(X)-110}.$$
One can easily find a constant $K_4'$ such that $g(C)\le 217$ provided $p_g(X)\ge K_4'$.
This contradicts to the assumption $g(C)\ge 648$. Thus, when $p_g(X)\ge K_4:=
\text{max}\{K_3,K_4'\}$,
$g(C)\le 647$.
\end{proof}

Theorems 1.5, 2.3 and 2.8 imply Theorem \ref{T:0.1}.

\begin{exmp}\label{E:2.9}
We give a family of examples where $X$ is canonically fibred by curves of genus $3$ and
the canonical map is of constant moduli. The original idea comes from \cite{Be1}. Let $W$
be a smooth projective surface with $p_g(W)=0$ and $C_0$ be a smooth curve of genus $2$.
Denote $X_0:=W\times C_0$. $p_1:X_0\longrightarrow W$ and
$p_2:X_0\longrightarrow C_0$ are projections. Suppose we have a divisor $H$ on $W$ such that $|2H|$ is base point free, $h^0(W, 2H)\ge 2$
and that $|K_W+H|$ gives a generically finite map.
Let $\theta=P-Q$ be a divisor on $C_0$ such that $2\theta\sim 0$. On $X_0$, denote $\delta:=p_1^*(H)+p_2^*(\theta)$. Then we can see that $\delta$ determines a smooth double cover $\pi:X\longrightarrow X_0$ with
$$K_X=\pi^*(K_{X_0}+\delta)$$
$$p_g(X)=h^0(X_0, p_1^*(K_W+H)).$$
Because $h^0(C_0, K_{C_0}+\theta)=1$, we can see that $\Phi_{K_X}$ factors through $\pi$ and $p_1$. If $K_W+H$ is nef, then $X$ is minimal. Denote $f:=p_1\circ\pi$. Then $f$ is the derived fibration from $\phi_1$ and a general fiber of $f$ is a curve of genus $3$.
\end{exmp}

\begin{setup} {\bf Open problem.}
It is quite interesting to consider a parallel problem. Let $X$ be a Gorenstein minimal projective
3-fold with at worst locally factorial terminal singularities. Suppose $\phi_1$ is generically finite onto its image. Is the generic
degree of $\phi_1$ universally upper bounded? We don't know any work on this problem.
It is well-known that the upper bound is 16 in surface case.
\end{setup}

\end{document}